\documentclass{siamart190516}

\usepackage{lipsum}
\usepackage{amsfonts}
\usepackage{graphicx}
\usepackage{epstopdf}
\usepackage{algorithmic}
\usepackage{amssymb}
\ifpdf
  \DeclareGraphicsExtensions{.eps,.pdf,.png,.jpg}
\else
  \DeclareGraphicsExtensions{.eps}
\fi

\newsiamremark{remark}{Remark}
\newsiamremark{hypothesis}{Hypothesis}
\crefname{hypothesis}{Hypothesis}{Hypotheses}
\newsiamthm{claim}{Claim}

\headers{Novel second-order fractional numerical formulas}{BaoLi Yin, Yang Liu, Hong Li and ZhiMin Zhang}

\title{Two families of novel second-order fractional numerical formulas and their applications to fractional differential equations\thanks{Submitted to the editors DATE.
\funding{The work of the second author was supported in part by the NSFC grant 11661058.
The work of the third author was supported in part by the NSFC grant 11761053, the NSF of Inner Mongolia 2017MS0107, and the program for Young Talents of Science and Technology in Universities of Inner Mongolia Autonomous Region NJYT-17-A07.
The work of the fourth author was supported in part by NSFC 11871092 and NSAF U1530401.
}}}

\author{
BaoLi Yin\thanks{School of Mathematical Sciences, Inner Mongolia University, Hohhot 010021,
China (\email{baolimath@aliyun.com}, \email{mathliuyang@imu.edu.cn}, \email{smslhong@imu.edu.cn}).}
\and Yang Liu\footnotemark[2]
\and Hong Li\footnotemark[2] 
\and ZhiMin Zhang\thanks{Beijing Computational Science Research Center, Beijing 100193, China (\email{zmzhang@csrc.ac.cn}); Department of Mathematics, Wayne State University, Detroit, MI 48202, USA (\email{zzhang@math.wayne.edu}).}
}

\usepackage{amsopn}

\ifpdf
\hypersetup{
  pdftitle={Two families of novel second-order fractional numerical formulas and their applications to fractional differential equations},
  pdfauthor={B.L. Yin, Y. Liu, Z.M. Zhang, and H. Li}
}
\fi

\begin{document}

\maketitle

\begin{abstract}
  In this article, we introduce two families of novel fractional $\theta$-methods by constructing some new generating functions to discretize the Riemann-Liouville fractional calculus operator $\mathit{I}^{\alpha}$ with a second order convergence rate.
A new fractional BT-$\theta$ method connects the fractional BDF2 (when $\theta=0$) with fractional trapezoidal rule (when $\theta=1/2$), and another novel fractional BN-$\theta$ method joins the fractional BDF2 (when $\theta=0$) with the second order fractional Newton-Gregory formula (when $\theta=1/2$).
To deal with the initial singularity, correction terms are added to achieve an optimal convergence order.
In addition, stability regions of different $\theta$-methods when applied to the Abel equations of the second kind are depicted, which demonstrate the fact that the fractional $\theta$-methods are A($\vartheta$)-stable.
Finally, numerical experiments are implemented to verify our theoretical result on the convergence analysis.
\end{abstract}



\section{Introduction}
\quad~ Fractional calculus is now an area attracting more and more attention both for its theory analysis interests and widespread applications in science and engineering fields. Many fractional derivatives such as Caputo type, Riemann-Liouville type, Riesz type lead to different fractional differential equations.
Considering the difficulties when solving equations with fractional calculus or the complex expressions of the analytic solutions, several popular numerical methods have been devised to efficiently get the numerical solutions. To formulate the numerical scheme of solving fractional differential equations, one need
to devise some efficient numerical formulas for fractional calculus operators.
Up to now, some high-order numerical approximations for fractional calculus operators, which have attracted a lot of attention, have been developed by some scholars; see
fractional linear multistep methods \cite{Lubich1,Jinbt1,Zeng2,Gaogh1,JinLiZhou},
L2-1${}_{\sigma}$ formula \cite{Alikhanov,Liaohl1}, WSGD operators \cite{Dengwh2,Zeng1,Liuydyw,Liuyd,Vong1}, and other high-order numerical schemes \cite{McLean1,McLean2,Xucj1,Ford1,Licp2,Zhengml1,Zhangzm1,Meerschaert1,Wangh1}. Here, we will consider some new second-order approximation formulas for Riemann-Liouville fractional calculus operators.
\par
First we state some definitions of the fractional calculus operators used in this paper.
The Riemann-Liouville fractional integral operator $\mathit{I}^{\alpha}$ is defined as
\begin{equation}\label{Intr.1}\begin{split}
\mathit{I}^{\alpha}u(x)=\frac{1}{\Gamma(\alpha)}\int_{0}^{x}(x-s)^{\alpha-1}u(s)\mathrm{d}s, \quad \text{for } \alpha>0,
\end{split}
\end{equation}
and set $\mathit{I}^0=I$, the identity operator.
The Riemann-Liouville fractional differential operator ${}^{RL}D_0^{\alpha}$ or  $\mathit{I}^{-\alpha} (\alpha>0)$ is defined as
\begin{equation}\label{Intr.2}\begin{split}
{}^{RL}D_0^{\alpha}u(x)
=\frac{\mathrm{d}^n}{\mathrm{d}x^n} \mathit{I}^{n-\alpha}u(x)
=\frac{1}{\Gamma(n-\alpha)}\frac{\mathrm{d}^n}{\mathrm{d}x^n}\int_{0}^{x}(x-s)^{n-\alpha-1}u(s)\mathrm{d}s,
\end{split}
\end{equation}
where $n=\lceil \alpha \rceil$. One can easily check that the Riemann-Liouville fractional calculus operators coincide with the classical ones when $\alpha$ takes integers.
\par
In \cite{Lubich1}, Lubich generalized the Dahlquist's convergence theorem for linear multistep methods to differential equations with fractional integral operator (\ref{Intr.1}) or fractional differential operator (\ref{Intr.2}). Three second-order numerical schemes, the fractional BDF2 (FBDF2), the fractional trapezoidal rule (FTR), and the generalized Newton-Gregory formula (GNGF2) have been devised with correction terms. Their corresponding generating functions are:
\begin{equation}\label{Intr.0.1}\begin{split}
\text{FBDF2}:&~\omega(\xi)=(3/2-2\xi+1/2 \xi^2)^{-\alpha},
\\
\text{FTR}:&~\omega(\xi)=2^{-\alpha}\big(\frac{1+\xi}{1-\xi}\big)^{\alpha},
\\
\text{GNGF2}:&~\omega(\xi)=\frac{1-\alpha/2+1/2 \alpha \xi}{(1-\xi)^{\alpha}}.
\end{split}
\end{equation}
It is natural to ask the connection between these three schemes based on (\ref{Intr.0.1}). In an early work, Liu et al. \cite{Liu1} proposed a BDF2-$\theta$ scheme based on these works \cite{Gaogh1,Wangyj1,Sunh} which connects the BDF2 and Crank-Nicolson (CN) scheme. Specifically, for a Cauchy problem: $y'=f(x,y)$ with initial condition $y(x_0)=y_0$, the BDF2-$\theta$ approximation formula is
\begin{equation}\label{Intr.0.2}\begin{split}
\frac{1}{2h}[(3-2\theta)u^n-(4-4\theta)u^{n-1}+(1-2\theta)u^{n-2}]=(1-\theta)f^n+\theta f^{n-1},
\end{split}
\end{equation}
where, $u^n$ is the approximation of $y(x_n)$ and $f^n=f(x_n,u^n)$.
One can easily check that condition $\theta=0$ implies the BDF2 and $\theta=\frac{1}{2}$ recovers the CN scheme (also known as the second-order Adams-Moulton method).
Based on this idea above and the fractional linear multistep methods developed by Lubich \cite{Lubich1}, we propose a family of novel fractional BT-$\theta$ method which connects the FBDF2 and FTR, and furthermore, we devise another family of approximation formula called new fractional BN-$\theta$ method which connects the FBDF2 and GNGF2.
There are at least two advantages for our novel fractional $\theta$ formulas.
(i) Since the convolution weights in our formulas depend on the parameter $\theta$, we can devise some formulas with special choices of this parameter to meet the assumptions of some techniques developed in literature for stability analysis of schemes.
(ii) From the aspect of numerical applications, we find that the FTR is superior to FBDF2 with a smaller error estimate and a better empirical convergence rate.
However, for fractional derivatives the FTR is not theoretically stable, see \cite{Lubich1}.
Now with the fractional BT-$\theta$ method, the bridge of the above two methods, we can practically take $\theta$ that is close to $\frac{1}{2}$ to obtain almost the best empirical results.

\par
Our main contributions are as follows:
\\
$\bigstar_1$ Propose a family of new fractional BT-$\theta$ method which generalizes the popular FBDF2 and FTR.
In addition, we devise another family of new fractional BN-$\theta$ method which connects the FBDF2 with GNGF2.
\\
$\bigstar_2$ Prove the convergence in detail by taking different techniques for these two families of novel fractional $\theta$ approximation formulas, discuss a correction technique for the problem with weak regularity solutions, and depict the stability regions of the proposed novel fractional $\theta$ formulas.
\\
$\bigstar_3$ Verify the convergence theories by choosing two numerical examples with smooth solutions and weak regularity solutions, respectively.
\par
The rest paper is outlined as follows:
In \cref{sec:two_methods}, we introduce the new fractional BT-$\theta$ method and the novel fractional BN-$\theta$ method with specific weights formulas and corresponding generating functions. A special case for the novel fractional BN-$\theta$ method when $\theta=1$ is listed out with a simple generating function which is different from all the three well known schemes (\ref{Intr.0.1}).
In \cref{sec:convergence}, we analyse the convergence of both fractional $\theta$-methods under the framework of Lubich \cite{Lubich1}. As one can see the analysis for the novel fractional BT-$\theta$ method is much easier than that of the other, since the former is derived directly by the linear multistep method (\ref{Intr.0.2}) while the latter is not.
In \cref{sec:stability}, we mainly analyse and depict the stability regions of the two methods when applied to a linear Abel integral equation of the second kind.
We also conduct some numerical experiments to confirm our theoretical analysis in \cref{sec:tests} with smooth and weak regularity solutions.
Finally, in \cref{sec:conc} we make some concluding remarks about the two families of novel fractional $\theta$-methods.
\section{Two families of new fractional $\theta$ methods}\label{sec:two_methods}
To derive the numerical scheme of the Riemann-Liouville fractional calculus operators $\mathit{I}^{\alpha}$, we first divide the interval $[0,L]$ into a uniform partition $0=x_0<x_1<\cdots<x_N=L$, with $h=L/N$ and $x_k=kh$ for $k=0,1,\cdots,N$.
For a given series $\{a_k\}_{k=0}^{\infty}$, we define the corresponding generating function $a(\xi)=\sum_{k=0}^{\infty}a_k \xi^k$.
And for a given generating function, we can also obtain the corresponding series $\{a_k\}_{k=0}^{\infty}$. Denote $u(x_n)$ by $u^n$ for simplicity.
\par
We need some definitions which can be found in \cite{Lubich1,Weilbeer} to clearly describe the fractional-$\theta$ methods.
For $\alpha \in \mathbb{R}$, we call the numerical approximation of $\mathit{I}^{\alpha} u$ at node $x_n$ given by
\begin{equation}\label{thet.1}
\begin{split}
{}_h\mathit{I}^{\alpha}u^n = h^{\alpha}\sum_{j=0}^{n}\omega_{n-j}u^j
+h^{\alpha}\sum_{j=1}^{s}\omega_{n,j}u^j,
\end{split}
\end{equation}
the fractional convolution quadrature $\omega$, where weights $\omega_j$ are called convolution weights and $\omega_{n,j}$ are called starting weights.
Denote ${}_h\Omega^{\alpha}u^n:=h^{\alpha}\sum_{j=0}^{n}\omega_{n-j}u^j$ as the convolution part.
In addition, define the convolution error $E^n:={}_h\Omega^{\alpha}u^n-\mathit{I}^{\alpha}u(x_n)$.
The formula for deriving the starting weights $\omega_{n,j}$ is stated in (\ref{conv.13}).
Note that for a sufficiently smooth solution, the starting part $h^{\alpha}\sum_{j=1}^{s}\omega_{n,j}u^j$ can be omitted. Nonetheless, an equation with fractional calculus operators shows the initial singularity; see\cite{Stynes1} and references therein.
\par
In what follows, based on (\ref{thet.1}), two families of novel second-order fractional $\theta$ methods by taking new convolution weights $\omega_{j}$ are proposed.
\subsection{New fractional BT-$\theta$ method}
Now we construct the first family of new fractional $\theta$ method by choosing the coefficients or convolution weights $\{\omega_j\}_{j=0}^{n}$ as follows
\begin{equation}\label{thet.2}
\begin{split}
\omega_j=&\bigg[\frac{2(1-\theta)}{3-2\theta}\bigg]^{\alpha}\sum_{k=0}^{j}\sum_{s=0}^{j-k}\kappa_k^{(1)}\kappa_s^{(2)}\kappa_{j-k-s}^{(3)},
\\
\kappa_k^{(1)}=&\bigg(\frac{\theta}{1-\theta}\bigg)^k
\begin{pmatrix}
  \alpha \\
  k
\end{pmatrix},
\kappa_k^{(2)}=\bigg(-\frac{1-2\theta}{3-2\theta}\bigg)^k
\begin{pmatrix}
  -\alpha \\
  k
\end{pmatrix},
\kappa_k^{(3)}=(-1)^k
\begin{pmatrix}
  -\alpha \\
  k
\end{pmatrix}.
\end{split}
\end{equation}
The parameter $\theta$ in (\ref{thet.2}) satisfies
\begin{equation}\label{thet.3}
\begin{split}
\theta \in (-\infty,\frac{1}{2}) \text{ when } \alpha \leq 0, \text{ and } \theta \in (-\infty,\frac{1}{2}] \text{ when } \alpha >0.
\end{split}
\end{equation}
One may check that the generating function $\omega(\xi)$ defined by (\ref{thet.2}) is
\begin{equation}\label{thet.4}
\begin{split}
\omega(\xi)=\bigg[\frac{1-\theta+\theta \xi}{(3/2-\theta)-(2-2\theta)\xi+(1/2-\theta)\xi^2}\bigg]^{\alpha}.
\end{split}
\end{equation}
Hence, when $\theta=0$, $\omega(\xi)=(\frac{3}{2}-2\xi+\frac{1}{2}\xi^2)^{-\alpha}$ is the FBDF2 and when $\theta=\frac{1}{2}$, $\omega(\xi)=2^{-\alpha}\big(\frac{1+\xi}{1-\xi}\big)^{\alpha}$ becomes FTR. Considering this aspect, we call the fractional convolution quadrature $\omega$ the fractional BT-$\theta$ method.
\subsection{New fractional BN-$\theta$ method}
In this subsection, we present another family of new second-order $\theta$ scheme, where the convolution weights $\{\omega_j\}_{j=0}^n$ are taken as
\begin{equation}\label{thet.5}
\begin{split}
\omega_j&=\bigg(\frac{2}{3-2\theta}\bigg)^{\alpha}\big[(1-\alpha \theta)\kappa_j+\alpha \theta \kappa_{j-1}\big],
\\
\kappa_j&=(-1)^j\sum_{k=0}^{j}\bigg(\frac{1-2\theta}{3-2\theta}\bigg)^{j-k}
\begin{pmatrix}
  -\alpha \\
  k
\end{pmatrix}
\begin{pmatrix}
  -\alpha \\
  j-k
\end{pmatrix},
\end{split}
\end{equation}
with $\theta \in (-\infty,1]$, $\alpha\theta \leq \frac{1}{2}$, and the generating function
\begin{equation}\label{thet.6}
\begin{split}
\omega(\xi)=\frac{1-\alpha\theta+\alpha\theta \xi}{\big[(3/2-\theta)-(2-2\theta)\xi+(1/2-\theta)\xi^2\big]^{\alpha}}.
\end{split}
\end{equation}
Now it is easy to check that for $\theta=0$, (\ref{thet.6}) reduces to FBDF2 and for $\theta=\frac{1}{2}$, (\ref{thet.6}) coincides with GNGF2. So, we define (\ref{thet.1}) with (\ref{thet.5}) as the fractional BN-$\theta$ method.
\par
 In particular, when $\theta=1$, we obtain the following generating function
\begin{equation}\label{thet.7}
\begin{split}
\omega(\xi)=2^{\alpha}\frac{1-\alpha+\alpha \xi}{(1-\xi^2)^{\alpha}},
\end{split}
\end{equation}
which leads to a new numerical scheme with a much simpler generating function.
\begin{remark}
Here, for the application of two families of novel fractional schemes, the convolution weights are provided in (\ref{thet.2}) and (\ref{thet.5}).
However the direct computation using the algorithms (\ref{thet.2}) and (\ref{thet.5}) will leads to bigger CPU time, so we provide another computing method in  \cref{sec:appendix} by which the complexity for deriving the weights $\{\omega_j\}_{j=0}^N$ is merely of $O(N)$.
\end{remark}

\section{Convergence analysis}\label{sec:convergence}
In this section, we analyse the convergence of the fractional convolution quadrature $\omega$ with weights defined in (\ref{thet.2}) and (\ref{thet.5}).
\par
First, we give the definition of convergence (to $\mathit{I}^{\alpha}$) of the fractional convolution quadrature $\omega$ of the form:
\begin{equation}\label{conv.0.1}
\omega(\xi)=r_1(\xi)^{\alpha}r_2(\xi)
\end{equation}
where $r_i(\xi)$ are rational functions.
As one can see, the generating functions (\ref{thet.4}) and (\ref{thet.6}) are special cases of (\ref{conv.0.1}).
\begin{definition}\label{Conv.defn.1}(Convergence, see \cite{Lubich1})
The fractional convolution quadrature $\omega$ with weights generated by (\ref{conv.0.1}) is convergent (to $\mathit{I}^{\alpha}$) of order $2$ if for all $\beta \in \mathbb{C}\setminus\{-1,-2,\cdots\}$,
\begin{equation}\label{conv.1}
  (\mathit{I}^{\alpha}x^{\beta})(x_n)-h^{\alpha}\sum_{j=0}^{n}\omega_{n-j}x_j^{\beta}
  =O(x_n^{\alpha+\beta-2}h^2)+O(x_n^{\alpha-1}h^{\beta+1}).
\end{equation}
\end{definition}
\subsection{Convergence for fractional BT-$\theta$ method}
We first examine the linear multistep method (LMM) in (\ref{Intr.0.2}),
which can be stated as, with $\theta \in (-\infty,1]$,
\begin{equation}\label{conv.2}
\begin{split}
u^n=\frac{4-4\theta}{3-2\theta}u^{n-1}
-\frac{1-2\theta}{3-2\theta}u^{n-2}
+h\bigg[\frac{2\theta}{3-2\theta}f^{n-1}+\frac{2(1-\theta)}{3-2\theta}f^n\bigg].
\end{split}
\end{equation}
The first and second characteristic polynomials of (\ref{conv.2}) are as the following
\begin{equation}\label{conv.3}
\begin{split}
\rho(r)=r^2-\frac{4-4\theta}{3-2\theta}r+\frac{1-2\theta}{3-2\theta}, \quad
\sigma(r)=\frac{2(1-\theta)}{3-2\theta}r+\frac{2\theta}{3-2\theta}.
\end{split}
\end{equation}
\begin{lemma}\label{Conv.lemma.2}
The LMM (\ref{conv.2}) is zero-stable and consistent of order $2$. Hence by Lax-Richtmyer theorem, (\ref{conv.2}) is convergent of order $2$ provided the error on the initial data tends to zeros as $O(h^2)$.
\end{lemma}
\textbf{Proof.} To prove (\ref{conv.2}) is consistent of order $2$, we just need to check the following equalities
\begin{equation}\label{conv.4}
\begin{split}
\begin{cases}
a_0+a_1=1,\\
-a_1+b_{-1}+b_0=1,\\
a_1+2b_{-1}=1,
\end{cases}
\end{split}
\end{equation}
where $a_0=\frac{4-4\theta}{3-2\theta}$, $a_1=-\frac{1-2\theta}{3-2\theta}$, $b_{-1}=\frac{2(1-\theta)}{3-2\theta}$, and $b_0=\frac{2\theta}{3-2\theta}$.
One can prove (\ref{conv.4}) by direct calculation.
\par
By the equivalence of zero-stability and root condition (see \cite{Quarteroni}, p.505), we need to check that
\begin{equation}\label{conv.5}
\begin{split}
\begin{cases}
|r_j|\leq 1,\quad j=0,1,\\
\rho'(r_j)\neq 0, \quad \text{if } |r_j|=1,
\end{cases}
\end{split}
\end{equation}
where $r_j$ are roots of $\rho(r)=0$. Indeed we have $r_0=1$ and $r_1=\frac{1-2\theta}{3-2\theta}$, and (\ref{conv.5}) holds for any $\theta \in (-\infty,1]$. The proof of the lemma is completed.
\par
We are now in a position to construct $\omega(\xi)$ in (\ref{thet.4}) by letting $\omega(\xi)=\big[\frac{\sigma(1/\xi)}{\rho(1/\xi)}\big]^{\alpha}$.
The following theorem shows that $\omega(\xi)$ is convergent of order $2$ for $\mathit{I}^{\alpha}$ provided condition (\ref{thet.3}) is satisfied.
\begin{theorem}\label{Conv.thm.1}
The fractional convolution quadrature $\omega$ with weights defined in (\ref{thet.2}) is convergent (to $\mathit{I}^{\alpha}$) of order $2$ provided condition (\ref{thet.3}) is satisfied.
\end{theorem}
\textbf{Proof.} With condition (\ref{thet.3}), we may easily prove that the root $\xi$ of $\sigma(\xi)=0$ satisfies $|\xi|<1$ when $\alpha \leq 0$, or $|\xi|\leq 1$ when $\alpha>0$.
Hence, by \cref{Conv.lemma.2} and Theorem 2.6 and example 2.9 on p.709 of \cite{Lubich1}, we have proved the fractional convolution quadrature $\omega$ with weights generated by (\ref{thet.4}) is convergent of order $2$ for $\mathit{I}^{\alpha}$. The proof for the theorem is completed.
\subsection{Convergence for fractional BN-$\theta$ method}
The following lemma reveals some facts about the order of magnitude of coefficients in functions having isolated singularity at initial value, which is crucial for the convergence analysis of the fractional BN-$\theta$ method.
\begin{lemma}\label{Conv.lem.3}(See \cite{Raisbeck})
If $-1<\alpha<1$, then
\begin{equation}\label{conv.5.1}
\begin{split}
\int_{0}^{\pi}x^{-\alpha} \cos nx \mathrm{d}x \sim n^{\alpha-1}\Gamma(1-\alpha)\cos\frac{1}{2}\pi(1-\alpha),
\end{split}
\end{equation}
and if $0<\alpha<2$, then
\begin{equation}\label{conv.5.2}
\begin{split}
\int_{0}^{\pi}x^{-\alpha} \sin nx \mathrm{d}x \sim n^{\alpha-1}\Gamma(1-\alpha)\sin\frac{1}{2}\pi(1-\alpha),
\end{split}
\end{equation}
where $\sim$ represents the ratio of the two sides approaches $1$ as $n$ approaches $\infty$.
\end{lemma}
\begin{theorem}\label{Conv.thm.2}
The fractional convolution quadrature $\omega$ with weights defined in (\ref{thet.5}) is convergent (to $\mathit{I}^{\alpha}$) of order $2$ under the condition $\theta \in (-\infty,1]$.
\end{theorem}
\textbf{Proof.} We will prove that the fractional convolution quadrature $\omega$ is stable and consistent of order $2$ for $\mathit{I}^{\alpha}$, then by Theorem 2.5 in \cite{Lubich1}, ${}_h \mathit{I}^{\alpha}$ is convergent (to $\mathit{I}^{\alpha}$) with second order convergence rate.
\par
By the relation $\omega_n=\big(\frac{2}{3-2\theta}\big)^{\alpha}\big[(1-\alpha \theta)\kappa_n+\alpha \theta \kappa_{n-1}\big]$, to prove (\ref{thet.6}) is stable (see Definition 2.1 in \cite{Lubich1}), i.e.,
\begin{equation}\label{conv.6}
\begin{split}
\omega_n=O(n^{\alpha-1}),
\end{split}
\end{equation}
we just need to prove that $\kappa_n=O(n^{\alpha-1})$, where $\{\kappa_n\}_{n=0}^{\infty}$ are the coefficients of
\begin{equation}\label{conv.7}
\begin{split}
\kappa(\xi)=&\big[(3/2-\theta)-(2-2\theta)\xi+(1/2-\theta)\xi^2\big]^{-\alpha}
\\=&(\frac{3}{2}-\theta)^{-\alpha}(1-\xi)^{-\alpha}(1-\gamma\xi)^{-\alpha},
\end{split}
\end{equation}
where $\gamma=\frac{1-2\theta}{3-2\theta} \in [-1,1)$, as $\theta \in (-\infty,1]$.
\par
Actually for any
$\alpha \in -\mathbb{N}= \{0,-1,-2,\cdots\}$, $\kappa_n=0$ for sufficiently large $n$ and (\ref{conv.6}) certainly is true.
Hence we mainly focus on $\alpha \in \mathbb{R}\setminus -\mathbb{N}$. Let $m=\lceil -\alpha \rceil$.
Now $\kappa(\xi)=(\frac{3}{2}-\theta)^{-\alpha}[(1-\xi)(1-\gamma \xi)]^m [(1-\xi)(1-\gamma \xi)]^{-\alpha-m}$, which means $\kappa_n=O(n^{\alpha-1})$ provided the coefficients of $[(1-\xi)(1-\gamma \xi)]^{-\alpha-m}$ are of $O(n^{\alpha+m-1})$. Noting that $m+\alpha \in (0,1)$, we next just analyse the case $\alpha \in (0,1)$ for (\ref{conv.7}).
\par
Based on the observation $\sum_{n=0}^{\infty}\kappa_n e^{-int}=\kappa(e^{-it})$, by fourier transform we have (see \cite{Podlubny})
\begin{equation}\label{conv.9}
\begin{split}
\kappa_n=\frac{1}{2\pi}\int_{0}^{2\pi}\kappa(e^{-it})e^{int}\mathrm{d}t.
\end{split}
\end{equation}
Note that $\kappa(e^{-it})$ is singular at $t=0$ and $t=2\pi$, hence we derive as
\begin{equation}\label{conv.10}
\begin{split}
\kappa_n=&\frac{1}{2\pi}\int_{0}^{\pi}\kappa(e^{-it})e^{int}\mathrm{d}t
+\frac{1}{2\pi}\int_{\pi}^{2\pi}\kappa(e^{-it})e^{int}\mathrm{d}t
\\=&\frac{1}{2\pi}\int_{0}^{\pi}\kappa(e^{-it})e^{int}
+\kappa(e^{it})e^{-int}\mathrm{d}t
\\=&\frac{1}{2\pi}\int_{0}^{\pi}\big[\kappa(e^{-it})+\kappa(e^{it})\big]\cos(nt)\mathrm{d}t
+\frac{i}{2\pi}\int_{0}^{\pi}\big[\kappa(e^{-it})-\kappa(e^{it})\big]\sin(nt)\mathrm{d}t.
\end{split}
\end{equation}
Careful calculations show that
\begin{equation}\label{conv.11}
\begin{split}
\kappa(e^{-it})+\kappa(e^{it}) \sim & 2\cos\frac{\alpha \pi}{2}(\frac{3}{2}-\theta)^{-\alpha}(1-\gamma)^{-\alpha}t^{-\alpha},
\\
i\big[\kappa(e^{-it})-\kappa(e^{it})\big] \sim & 2\sin\frac{\alpha \pi}{2}(\frac{3}{2}-\theta)^{-\alpha}(1-\gamma)^{-\alpha}t^{-\alpha}.
\end{split}
\end{equation}
Now with \cref{Conv.lem.3} we have proved the stability of the fractional convolution quadrature $\omega$ with weights generated by (\ref{thet.6}).
\par
By Lemma 3.2 in \cite{Lubich1}, consistency of (\ref{thet.6}) is equivalent to the condition that expansion of $\omega(\xi)=(1-\xi)^{-\alpha}[c_0+c_1(1-\xi)+(1-\xi)^2 \tilde{r}(\xi)]$ satisfies
$c_0=\gamma_0$ and $c_1=\gamma_1$, where $\gamma_i$ denotes the coefficients of
\begin{equation}\label{conv.11.1}
\begin{split}
\sum_{i=0}^{\infty} \gamma_i(1-\xi)^i=\bigg(-\frac{\ln \xi}{1-\xi}\bigg)^{-\alpha}.
\end{split}
\end{equation}
Direct calculation shows $\gamma_0=1$ and $\gamma_1=-\frac{\alpha}{2}$. For $\omega(\xi)$ in (\ref{thet.6}), we have $\omega(\xi)=(1-\xi)^{-\alpha}\tilde{\omega}(\xi)$ with $\tilde{\omega}(\xi)=(\frac{3}{2}-\theta)^{-\alpha}(1-\alpha \theta+\alpha \theta \xi)(1-\frac{1-2\theta}{3-2\theta} \xi)^{-\alpha}$. Now expanding $\tilde{\omega}(\xi)$ at $\xi=1$, we can easily derive that $c_0=\tilde{\omega}(1)=1$ and $c_1=-\tilde{\omega}'(1)=-\frac{\alpha}{2}$, which show the claim of consistency. The proof of the theorem is completed.
\begin{remark}\label{rem.1}
We remark that although condition $\alpha \theta \leq \frac{1}{2}$ is not used in the proof for fractional BN-$\theta$ method, it is not A($\frac{\pi}{2}$)-stable if $\alpha \theta > \frac{1}{2}$ under the condition $\alpha \in (0,1)$. Generally, we may choose $\theta$ such that $\alpha \theta \leq \frac{1}{2}$ and $\theta \in (-\infty,1]$ to assure A($\vartheta$)-stability. See section $4$ for more information.
\end{remark}
\subsection{Correction technique}
For a smooth function $u$, the convolution error $E^n:={}_h\Omega^{\alpha}u^n-\mathit{I}^{\alpha}u(x_n)$ of the two families of fractional $\theta$-methods is $O(h^2)$ uniformly for bounded $x$. However, for solutions with weak regularity, the convergence rate $O(h^2)$ cannot be maintained. To overcome this difficulty, we apply the technique of adding correction terms introduced by Lubich \cite{Lubich1} to our fractional $\theta$-methods.
\begin{lemma}\label{Conv.lemma.1}(See \cite{Lubich1})
If the fractional convolution quadrature $\omega$ with weights generated by (\ref{conv.0.1}) is convergent of order $2$ for $\mathit{I}^{\alpha}$, then, for every $\beta \in \mathbb{C}\setminus\{-1,-2,\cdots\}$, there exist starting weights $\omega_{n,j}$ such that
\begin{equation}\label{conv.12}
\begin{split}
\mathit{I}^{\alpha}u \big|_{x=x_n}
= h^{\alpha}\sum_{j=0}^{n}\omega_{n-j}u^j
+h^{\alpha}\sum_{j=1}^{s}\omega_{n,j}u^j
+O(h^2),
\end{split}
\end{equation}
for any function $u(x)=x^{\beta}f(x)$ with $f$ sufficiently differentiable,
where the starting weights $\omega_{n,j}$ satisfy a linear system,
\begin{equation}\label{conv.13}
\sum_{j=1}^{s}\omega_{n,j} j^{\ell} = \frac{\Gamma(\ell+1)}{\Gamma(\ell+\alpha+1)}n^{\ell+\alpha}
-\sum_{j=0}^{n}\omega_{n-j}j^{\ell},
\end{equation}
with $\ell \in \Sigma=\{\ell=\beta+q:q\in \mathbb{N},\ell < 2-\min\{1,\alpha\}\}$ and $s:=\text{card }\Sigma$.
\end{lemma}
We have proved that the two families of fractional $\theta$-methods in this paper are convergent of order $2$, hence by \cref{Conv.lemma.1} the correction technique can be applied to our fractional $\theta$-methods.
\begin{remark}
The asymptotic properties of $\omega_n$ of the two families of fraction $\theta$-methods and the starting weights $\omega_{n,j}$ derived by (\ref{conv.13}) are
\begin{equation}\label{conv.14}
\omega_n=O(n^{\alpha-1}), \quad \omega_{n,j}=O(n^{\alpha-\gamma}) \text{ with } \gamma=\min\{1,\alpha\},
\end{equation}
respectively.
\end{remark}

\section{Stability region}\label{sec:stability}
In this section, we will apply the two families of fractional $\theta$-methods to a linear Abel integral equation of the second kind to explore the stability region with different parameter $\theta$. The model equation is:
\begin{equation}\label{stab.1}
\begin{split}
u(x)=f(x)+\frac{\lambda}{\Gamma(\alpha)}\int_{0}^{x}(x-s)^{\alpha-1}u(s)\mathrm{d}s, \quad 0 < \alpha <1.
\end{split}
\end{equation}
The analytic stability region for (\ref{stab.1}) is $|\arg \lambda -\pi|<(1-\frac{1}{2}\alpha)\pi$, i.e., for any $\lambda \in \mathbb{C}$ within the region, the analytic solution $u(x)$ satisfies $u(x)\to 0$ as $x\to \infty$ whenever $f(x)$ converges to a finite limit (see \cite{Lubich2}).
Applying the fractional $\theta$-methods to (\ref{stab.1}) and assuming $u$ is sufficiently smooth, we have
\begin{equation}\label{stab.2}
\begin{split}
u^n=f(x_n)+\lambda h^{\alpha}\sum_{j=0}^{n}\omega_{n-j}u^j.
\end{split}
\end{equation}
\par
First we state some definitions concerning the stability of equation (\ref{stab.1}).
\begin{definition}(See \cite{Lubich2})\label{defn.1}
Assume sequence $\{f(x_n)\}_{n=0}^{\infty}$ has a finite limit. A fractional convolution quadrature $\omega$ for the Abel equation (\ref{stab.1}) is called A-stable if the numerical solution $u^n$ given by (\ref{stab.2}) satisfies
\begin{equation}\label{stab.2.0}
u^n \to 0 \quad \text{as} \quad n\to \infty
\end{equation}
for any $h>0$ and any $\lambda$ in the analytical stability region $|\arg \lambda-\pi|<(1-\frac{1}{2}\alpha)\pi$.
\end{definition}
\begin{definition}(See \cite{Lubich2})\label{defn.2}
The stability region $S$ of a fractional convolution quadrature $\omega$ is the set of all complex $z=\lambda h^{\alpha}$ for which the numerical solution $u^n$ given by (\ref{stab.2}) satisfies
\begin{equation}\label{stab.3}
u^n \to 0 \quad as \quad n \to \infty
\end{equation}
provided $\{f(x_n)\}$ has a finite limit. Moreover, $\omega$ is called A($\vartheta$)-stable if $S$ contains the sector $|\arg z-\pi|<\vartheta$.
\end{definition}
\begin{theorem}\label{stab.thm.1}
For the fractional convolution quadrature $\omega$ with weights defined in (\ref{thet.2}) or (\ref{thet.5}), the stability region $S$ is
\begin{equation}\label{stab.4}
S=\mathbb{C}\setminus\{1/\omega(\xi):|\xi|\leq1\}.
\end{equation}
\end{theorem}
\textbf{Proof.} Considering the Theorem 2.1 in \cite{Lubich2}, a sufficient condition for (\ref{stab.4}) is to express $\omega_n$ as
\begin{equation}\label{stab.5}
\omega_n=(-1)^n
\begin{pmatrix}
  -\alpha \\
  n
\end{pmatrix}+v_n
\quad \text{with }\{v_n\} \in \ell^1.
\end{equation}
Actually we have $\omega(\xi)=(1-\xi)^{-\alpha}\tilde{\omega}(\xi)$, where
\begin{equation*}\label{stab.6}
\begin{split}
\tilde{\omega}(\xi)=\bigg[\frac{1-\theta+\theta\xi}{(\frac{1}{2}-\theta)(1-\xi)+1}\bigg]^{\alpha}&, \quad \text{for BT-$\theta$ method,}
\end{split}
\end{equation*}
and
\begin{equation*}\label{stab.6.1}
\begin{split}
\tilde{\omega}(\xi)=\bigg(\frac{3}{2}-\theta\bigg)^{-\alpha}(1-\alpha \theta+\alpha \theta \xi)\bigg(1-\frac{1-2\theta}{3-2\theta} \xi\bigg)^{-\alpha}&, \quad \text{for BN-$\theta$ method.}
\end{split}
\end{equation*}
Expanding $\tilde{\omega}(\xi)$ at $1$, we have
\begin{equation}\label{stab.7}
\begin{split}
\omega(\xi)=(1-\xi)^{-\alpha}\bigg[1+\frac{\alpha}{2}(1-\xi)\bigg]+(1-\xi)^2 r(\xi).
\end{split}
\end{equation}
By Lemma 3.3 in \cite{Lubich1}, the coefficients $r_n$ of $r(\xi)$ in (\ref{stab.7}) satisfy $r_n=O(n^{\alpha-1})$, since the fractional convolution quadrature $\omega$ is convergent. Now (\ref{stab.7}) means (\ref{stab.5}) is true and the proof of the theorem is completed.
\par
We next depict the stability regions with different $\theta$ and $\alpha$ for both of the methods. According to \cref{stab.thm.1}, we have the relation $1/\omega(-1) \in (\mathbb{C}\setminus S) \cap \mathbb{R}$, which means the set $\mathbb{C}\setminus S$ must at least contain the point $1/\omega(-1)$ (lie at the $x$-axis).
Direct calculations show that
\begin{equation}\label{stab.8}
1/\omega(-1)=
\begin{cases}
  4^{\alpha}(1-\theta)^{\alpha}/(1-2\theta)^{\alpha}, & \mbox{for BT-$\theta$ method,}\\
  4^{\alpha}(1-\theta)^{\alpha}/(1-2\alpha \theta), & \mbox{for BN-$\theta$ method}.
\end{cases}
\end{equation}
\begin{figure}[h]
\begin{center}
\begin{minipage}{6cm}
  \centering\includegraphics[width=6cm]{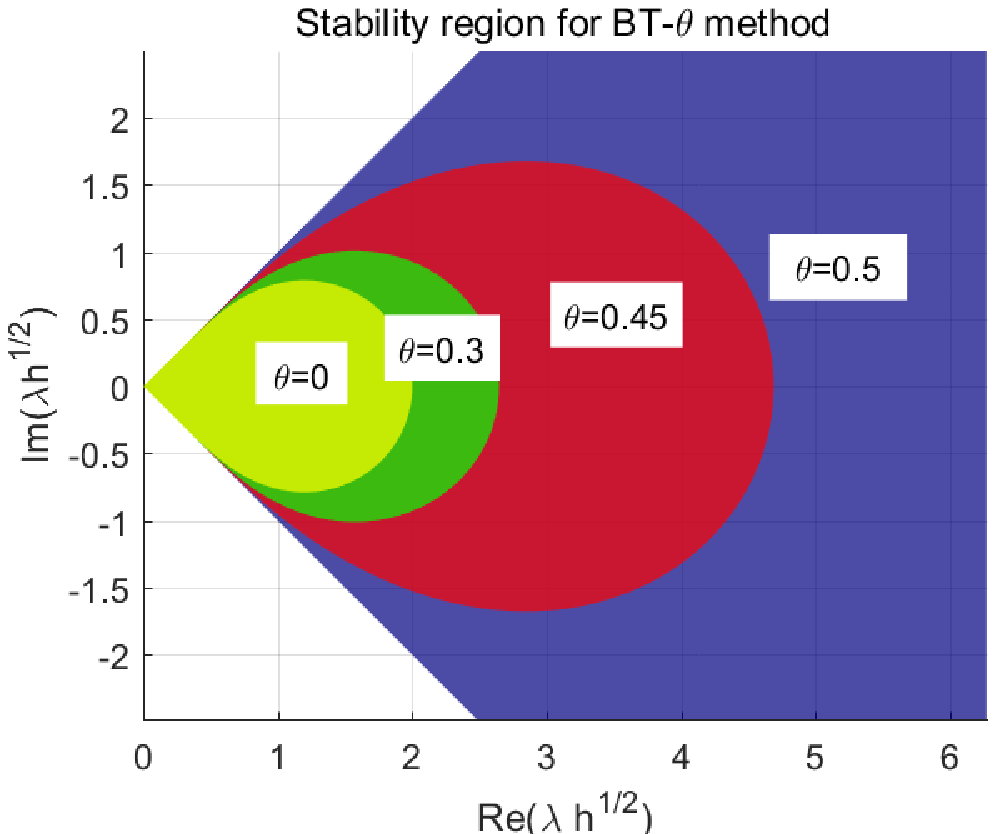}
  \caption{$\alpha=\frac{1}{2}$, $\theta \in [0,\frac{1}{2}]$.}\label{C1}
\end{minipage}
\begin{minipage}{6cm}
  \centering\includegraphics[width=6cm]{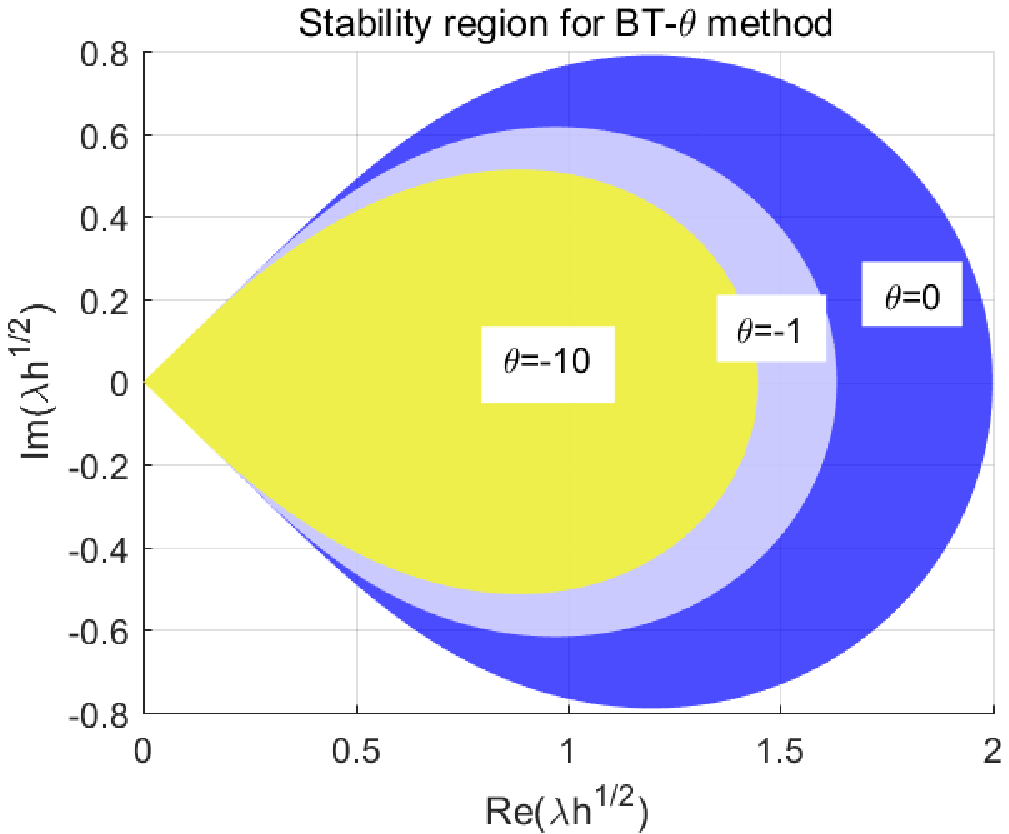}
  \caption{$\alpha=\frac{1}{2}$, $\theta \in (-\infty,0]$.}\label{C2}
\end{minipage}
\end{center}
\end{figure}
\begin{figure}[h]
\begin{center}
\begin{minipage}{6cm}
  \centering\includegraphics[width=6cm]{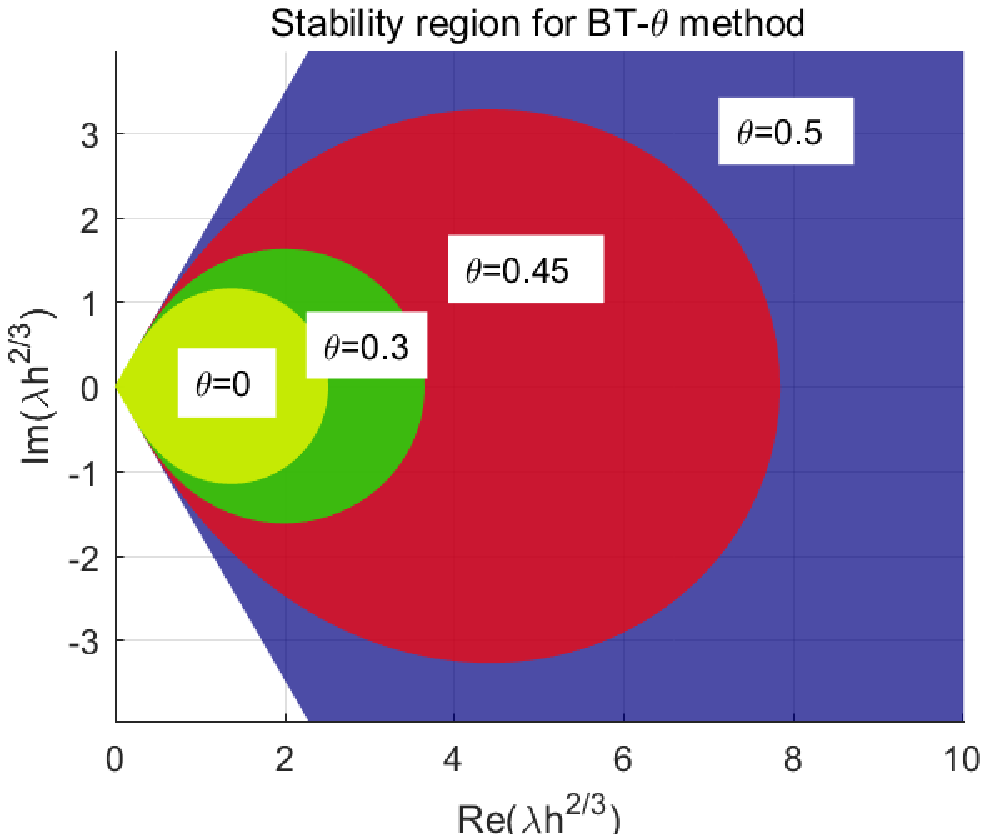}
  \caption{$\alpha=\frac{2}{3}$, $\theta \in [0,\frac{1}{2}]$.}\label{C3}
\end{minipage}
\begin{minipage}{6cm}
  \centering\includegraphics[width=6cm]{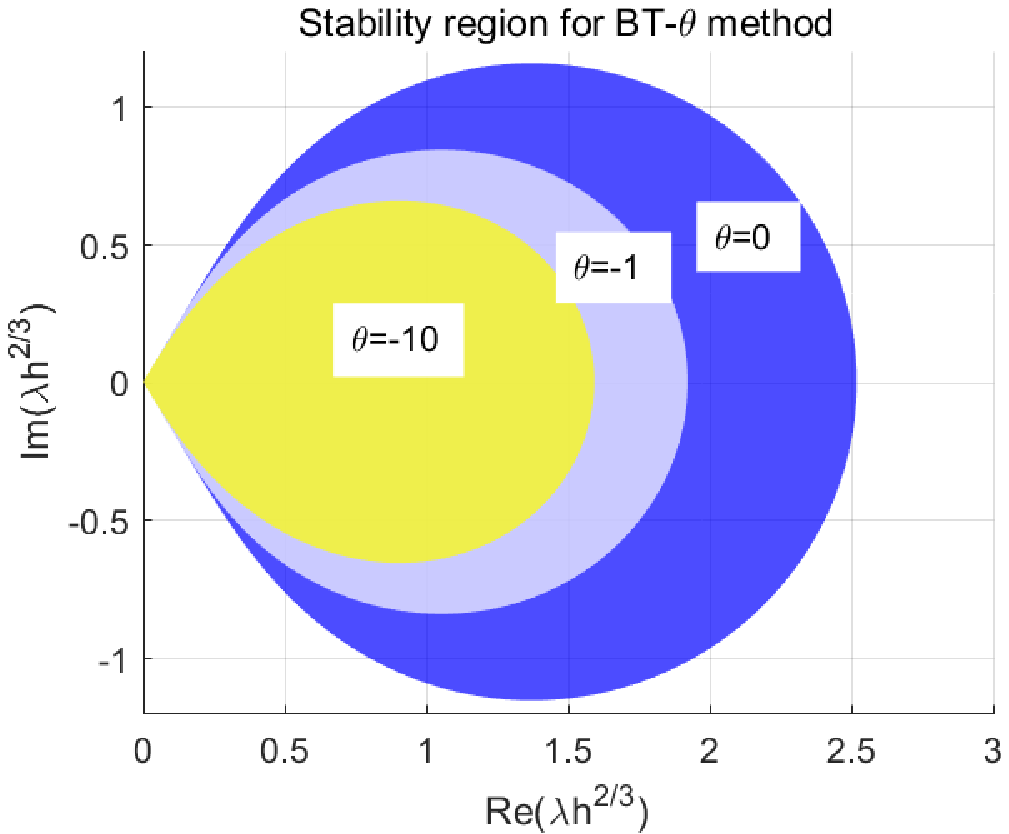}
  \caption{$\alpha=\frac{2}{3}$, $\theta \in (-\infty,0]$.}\label{C4}
\end{minipage}
\end{center}
\end{figure}
\subsection{Stability regions for the fractional BT-$\theta$ method}
For the fractional BT-$\theta$ method with $\theta \in (-\infty,\frac{1}{2}]$, it holds that $1/\omega(-1) \in (2^{\alpha},+\infty)$. Figs. \ref{C1}-\ref{C4} show that $S$ contains the analytical stability region $|\arg \lambda-\pi|<(1-\frac{1}{2}\alpha)\pi$
(In Figs. \ref{C1}-\ref{C2}, the analytical stability region is $|\arg \lambda-\pi|<\frac{3}{4}\pi$ and in Figs. \ref{C3}-\ref{C4}, the region is $|\arg \lambda-\pi|<\frac{2}{3}\pi$. The shaded area represents set $\mathbb{C}\setminus S$).
Hence, by the \cref{defn.1}, the fractional BT-$\theta$ method is A-stable for $\theta \in (-\infty,\frac{1}{2}]$.
One may also find out that $1/\omega(-1)$ is exactly one of the points where shaded area intersects with $x$-axis (also see Theorem 4.1 in \cite{Lubich2}).

\subsection{Stability regions for the fractional BN-$\theta$ method}
\begin{figure}[h]
\begin{center}
\begin{minipage}{6cm}
  \centering\includegraphics[width=6cm]{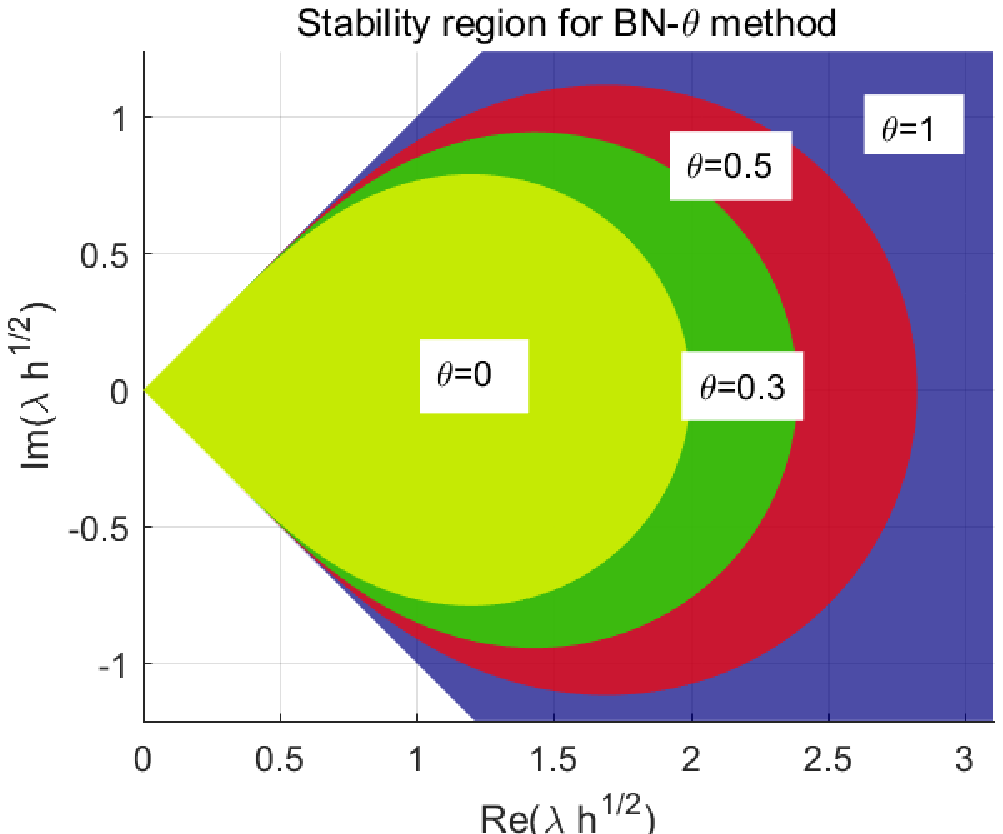}
  \caption{$\alpha=\frac{1}{2}$, $\theta \in [0,1]$.}\label{C5}
\end{minipage}
\begin{minipage}{6cm}
  \centering\includegraphics[width=6cm]{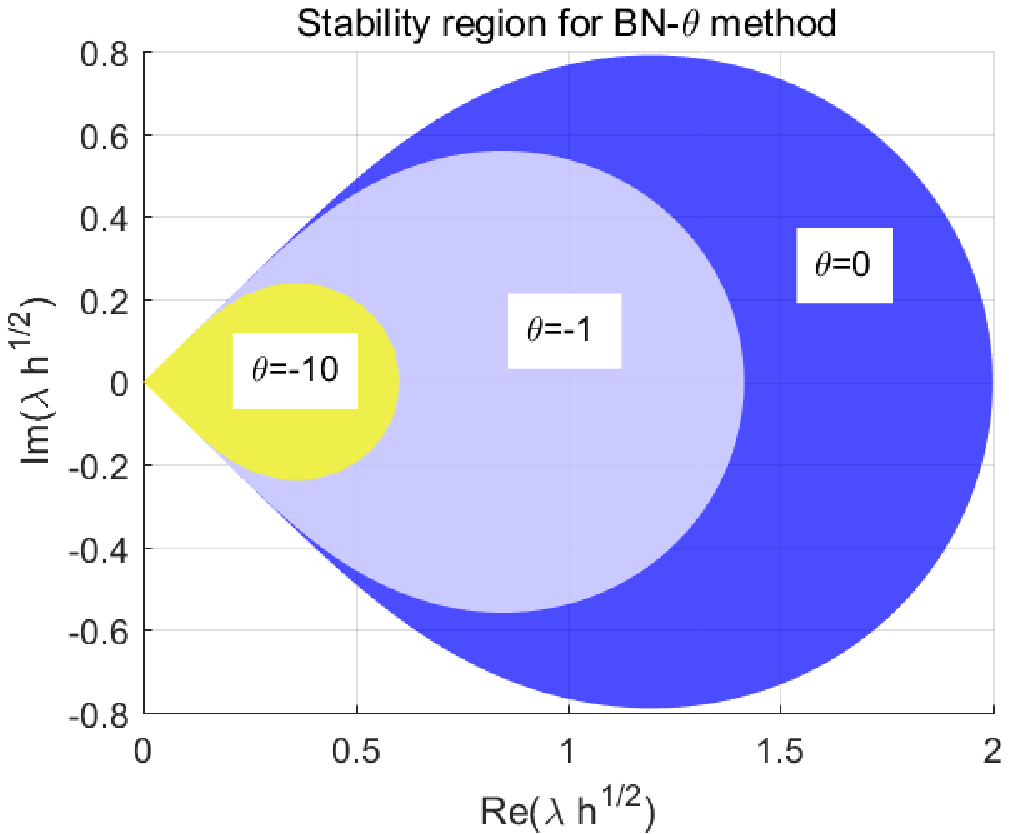}
  \caption{$\alpha=\frac{1}{2}$, $\theta \in (-\infty,0]$.}\label{C6}
\end{minipage}
\end{center}
\end{figure}

\begin{figure}[h]
\begin{center}
\begin{minipage}{6cm}
  \centering\includegraphics[width=6cm]{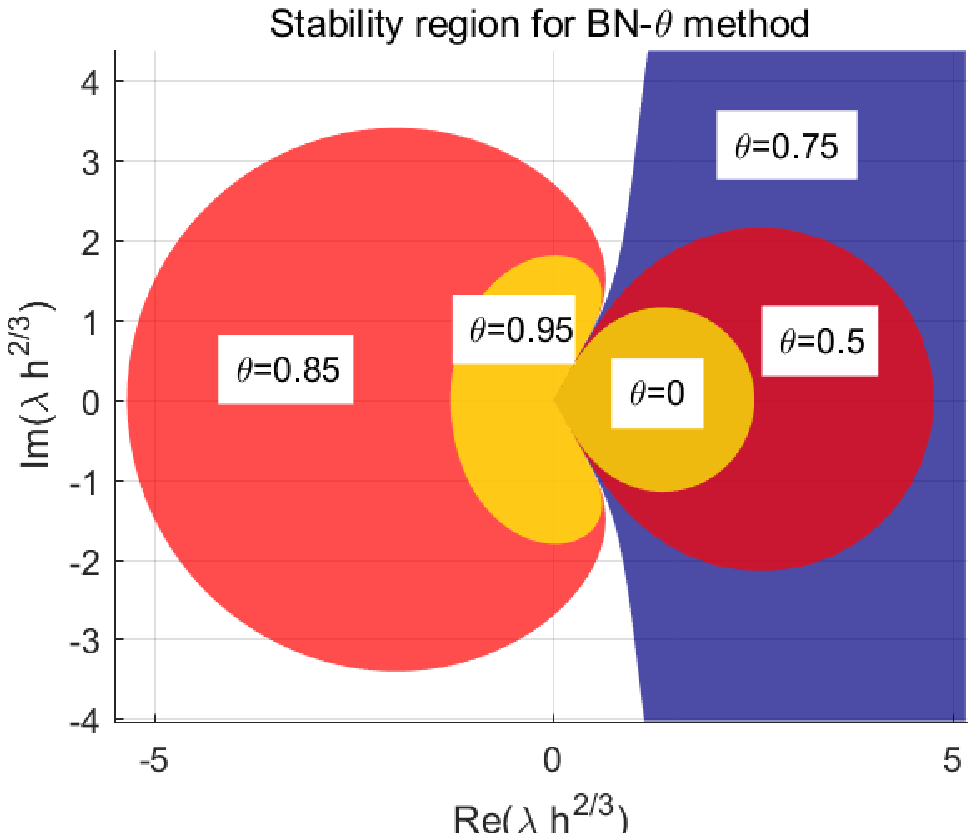}
  \caption{$\alpha=\frac{2}{3}$, $\theta \in [0,1]$.}\label{C7}
\end{minipage}
\begin{minipage}{6cm}
  \centering\includegraphics[width=6cm]{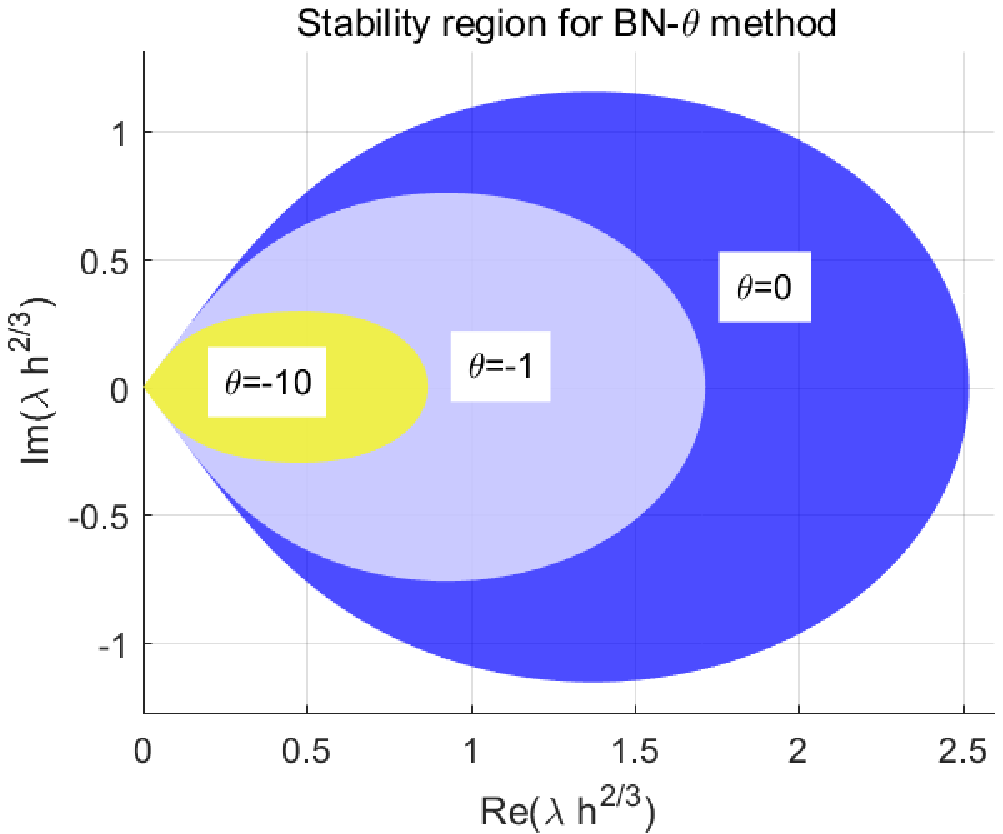}
  \caption{$\alpha=\frac{2}{3}$, $\theta \in (-\infty,0]$.}\label{C8}
\end{minipage}
\end{center}
\end{figure}

\begin{figure}[h]
\begin{center}
\begin{minipage}{6cm}
  \centering\includegraphics[width=6cm]{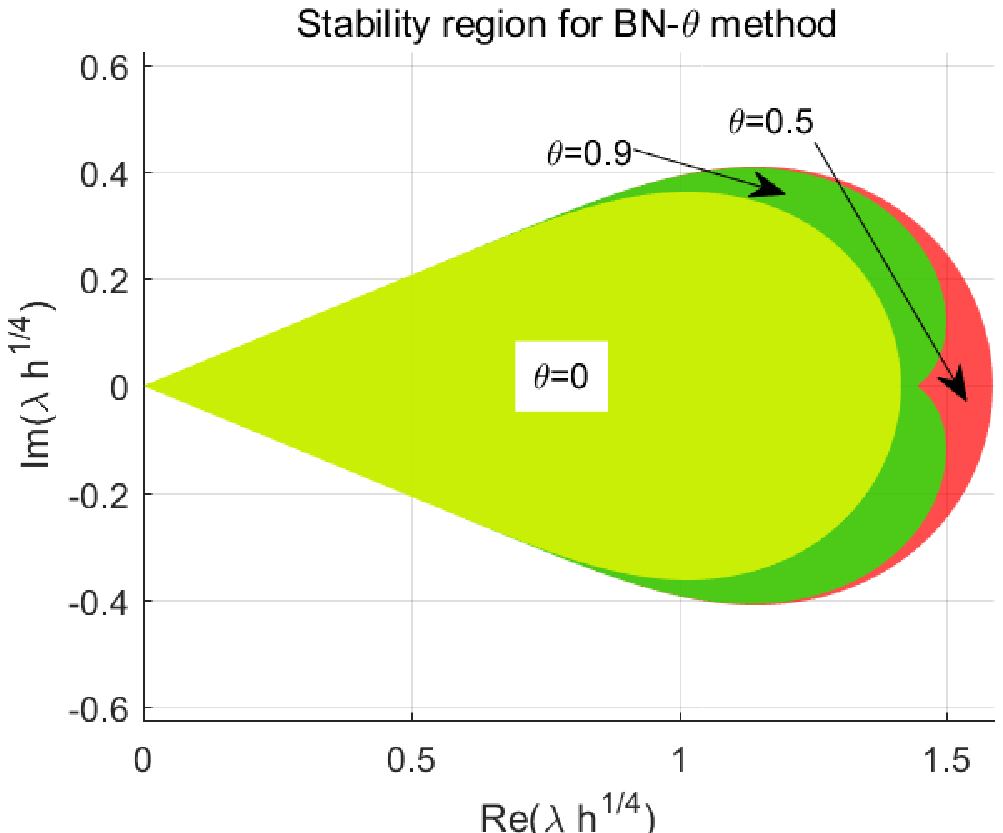}
  \caption{$\alpha=\frac{1}{4}$, $\theta \in [0,1]$.}\label{C9}
\end{minipage}
\begin{minipage}{6cm}
  \centering\includegraphics[width=6cm]{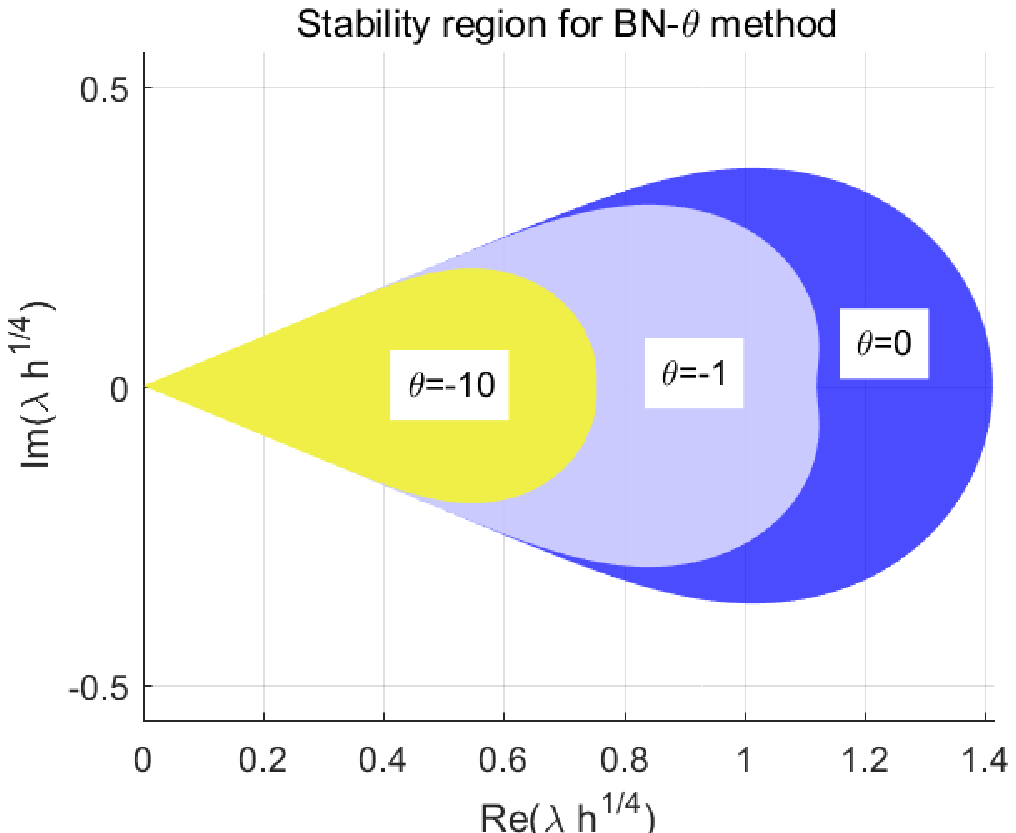}
  \caption{$\alpha=\frac{1}{4}$, $\theta \in (-\infty,0]$.}\label{C10}
\end{minipage}
\end{center}
\end{figure}

The stability region for this method is complicated as one can see from (\ref{stab.8}) that $1/\omega(-1)<0$ if $\alpha \theta>\frac{1}{2}$, in which case the method is not A($\frac{\pi}{2}$)-stable.
Fig. \ref{C5} and Fig. \ref{C6} show some similar properties of the stability region when compared with Fig. \ref{C1} and Fig. \ref{C2}, respectively.
However, the stability region for $\theta=-10$ in Fig. \ref{C6} is much smaller than the one in Fig. \ref{C2}.
When $\alpha > \frac{1}{2}$, we have chosen $\theta \in [0,1]$ in Fig. \ref{C7} such that $\alpha \theta>\frac{1}{2}$, in which case the method is not A($\frac{\pi}{2}$)-stable.
Fig. \ref{C8} also confirms the fact that for $\alpha \in (0,1)$ and $\theta \leq 0$, the fractional BN-$\theta$ method is A($\frac{\pi}{2}$)-stable.
When taking $\alpha=\frac{1}{4}$, Fig. \ref{C9} and Fig. \ref{C10} show some different properties of the shape of the stability regions, such as $\theta=0.9$ or $\theta=-1$.
\par
All in all, we may conclude from the above stability regions with different $\alpha$ and $\theta$, that fractional BT-$\theta$ method is A-stable for any $\theta \in (-\infty,\frac{1}{2}]$, and fractional BN-$\theta$ method is A($\frac{\pi}{2}$)-stable provided $\theta \in (-\infty,\min\{1,\frac{1}{2\alpha}\}]$ for the equation (\ref{stab.1}).

\section{Numerical tests}\label{sec:tests}
In this section, we take some numerical experiments to verify the efficiency of the proposed fractional $\theta$-methods. The first example assumes solution is sufficiently smooth in which case (\ref{thet.1}) is used, and the second example assumes some weak regularity on the solution in which case (\ref{conv.7}) is used with correction terms.
\subsection{Example with a sufficiently smooth solution}

We consider the following linear Caputo fractional ODE:
\begin{equation}\label{nume.1}
\begin{split}
{}^C D_{0}^{\alpha} u=u+f(x),
\quad x \in (0,L] \quad \text{and } u(0)=u_0,
\end{split}
\end{equation}
where $\alpha \in (0,1)$ and ${}^C D_{0}^{\alpha}$ is the Caputo differential operator.
Under the condition $\alpha \in (0,1)$, we have ${}^C D_{0}^{\alpha} u=\mathit{I}^{-\alpha}(u-u_0)$.
Let $v(x)=u(x)-u_0$.
Now (\ref{nume.1}) is formulated as the following Riemann-Liouville fractional ODE
\begin{equation}\label{nume.2}
\begin{split}
\mathit{I}^{-\alpha}v=v+f(x)+u_0,
\quad x \in (0,L] \quad \text{and } v(0)=0.
\end{split}
\end{equation}
\par
We take $u=1+x^3$, and the corresponding $f(x)$ is $f(x)=\frac{6x^{3-\alpha}}{\Gamma(4-\alpha)}-x^3-1$. By formula (\ref{thet.1}) we can easily derive the numerical scheme as
\begin{equation}\label{nume.3}
\begin{split}
{}_h\Omega^{-\alpha}v^n=v^n+f(x_n)+u_0.
\end{split}
\end{equation}
\par
In Table \ref{tab1}, we denote Error($\theta$) as $\max_{1 \leq n \leq N}|U^n-u(x_n)|$ where $U^n$ is the numerical solution obtained by the fractional BT-$\theta$ method. For different $\theta$ taken as $-1,0, 0.2$ or $0.45$ and different $\alpha \in (0,1)$, we have obtained a second-order convergence rate as desired. A phenomenon is that for a smaller $\theta$, we get a larger error, and hence the convergence rate may be also affected.
\par
In Table \ref{tab2}, we apply the fractional BN-$\theta$ method to (\ref{nume.1}) with choice of $\theta$  satisfying $\theta \leq 1$ and $-\alpha\theta \leq \frac{1}{2}$. One can see the convergence rate is also $O(h^2)$, which is in line with our theoretical result.
\begin{table}[tbhp]
{\footnotesize
\caption{The convergence rate for BT-$\theta$ method.}\label{tab1}
\begin{center}
\begin{tabular}{|c|l|c|c|c|c|c|c|c|c|} \hline
$\alpha$ &  $h$ & Error(-1) & rate &  Error(0) & rate &  Error(0.2) & rate &  Error(0.45) & rate \\ \hline
	&	   1/4  	&	2.168E-01	&	---     &	1.190E-01	&	---     &	9.048E-02	&	---     &	4.865E-02	&	--- \\
	&	   1/8  	&	6.772E-02	&	1.68 	&	3.309E-02	&	1.85 	&	2.430E-02	&	1.90 	&	1.225E-02	&	1.99 	\\
0.1	&	   1/16 	&	1.984E-02	&	1.77 	&	8.869E-03	&	1.90 	&	6.370E-03	&	1.93 	&	3.088E-03	&	1.99 	\\
	&	   1/32 	&	5.437E-03	&	1.87 	&	2.306E-03	&	1.94 	&	1.636E-03	&	1.96 	&	7.765E-04	&	1.99 	\\
	&	   1/64 	&	1.428E-03	&	1.93 	&	5.886E-04	&	1.97 	&	4.148E-04	&	1.98 	&	1.948E-04	&	2.00 	\\
\hline
	&	   1/4  	&	2.622E-01	&	---     &	1.419E-01	&	---     &	1.073E-01	&	---     &	5.713E-02	&	--- \\
	&	   1/8  	&	8.108E-02	&	1.69 	&	3.915E-02	&	1.86 	&	2.865E-02	&	1.90 	&	1.436E-02	&	1.99 	\\
0.5	&	   1/16 	&	2.352E-02	&	1.79 	&	1.044E-02	&	1.91 	&	7.482E-03	&	1.94 	&	3.617E-03	&	1.99 	\\
	&	   1/32 	&	6.405E-03	&	1.88 	&	2.705E-03	&	1.95 	&	1.917E-03	&	1.96 	&	9.087E-04	&	1.99 	\\
	&	   1/64 	&	1.676E-03	&	1.93 	&	6.894E-04	&	1.97 	&	4.856E-04	&	1.98 	&	2.278E-04	&	2.00 	\\
\hline
	&	   1/4  	&	3.111E-01	&	---     &	1.658E-01	&	---     &	1.247E-01	&	---     &	6.561E-02	&	--- \\
	&	   1/8  	&	9.544E-02	&	1.70 	&	4.557E-02	&	1.86 	&	3.325E-02	&	1.91 	&	1.656E-02	&	1.99 	\\
0.9	&	   1/16 	&	2.747E-02	&	1.80 	&	1.210E-02	&	1.91 	&	8.660E-03	&	1.94 	&	4.173E-03	&	1.99 	\\
	&	   1/32 	&	7.437E-03	&	1.88 	&	3.129E-03	&	1.95 	&	2.215E-03	&	1.97 	&	1.048E-03	&	1.99 	\\
	&	   1/64 	&	1.940E-03	&	1.94 	&	7.962E-04	&	1.97 	&	5.606E-04	&	1.98 	&	2.628E-04	&	2.00 	\\
\hline
\end{tabular}
\end{center}
}
\end{table}
\begin{table}[tbhp]
{\footnotesize
\caption{The convergence rate for BN-$\theta$ method.}\label{tab2}
\begin{center}
\begin{tabular}{|c|l|c|c|c|c|c|c|c|c|} \hline
$\alpha$ &  $h$ & Error(-0.5) & rate &  Error(0) & rate &  Error(0.5) & rate &  Error(1) & rate\\
\hline
	&	   1/4  	&	2.204E-01	&	---     &	1.190E-01	&	---     &	8.810E-02	&	---     &	1.487E-01	&	--- \\
	&	   1/8  	&	6.476E-02	&	1.77 	&	3.309E-02	&	1.85 	&	2.337E-02	&	1.91 	&	3.967E-02	&	1.91 	\\
0.1	&	   1/16 	&	1.813E-02	&	1.84 	&	8.869E-03	&	1.90 	&	6.081E-03	&	1.94 	&	1.041E-02	&	1.93 	\\
	&	   1/32 	&	4.840E-03	&	1.91 	&	2.306E-03	&	1.94 	&	1.555E-03	&	1.97 	&	2.679E-03	&	1.96 	\\
	&	   1/64 	&	1.253E-03	&	1.95 	&	5.886E-04	&	1.97 	&	3.935E-04	&	1.98 	&	6.804E-04	&	1.98 	\\
\hline
	&	   1/4  	&	2.903E-01	&	---     &	1.419E-01	&	---     &	1.220E-01	&	---     &	2.493E-01	&	--- \\
	&	   1/8  	&	8.326E-02	&	1.80 	&	3.915E-02	&	1.86 	&	3.287E-02	&	1.89 	&	6.786E-02	&	1.88 	\\
0.5	&	   1/16 	&	2.300E-02	&	1.86 	&	1.044E-02	&	1.91 	&	8.630E-03	&	1.93 	&	1.812E-02	&	1.91 	\\
	&	   1/32 	&	6.097E-03	&	1.92 	&	2.705E-03	&	1.95 	&	2.218E-03	&	1.96 	&	4.711E-03	&	1.94 	\\
	&	   1/64 	&	1.573E-03	&	1.95 	&	6.894E-04	&	1.97 	&	5.627E-04	&	1.98 	&	1.203E-03	&	1.97 	\\
\hline
	&	   1/4  	&	3.531E-01	&	---     &	1.598E-01	&	---     &	1.510E-01	&	---     &	3.361E-01	&	--- \\
	&	   1/8  	&	9.930E-02	&	1.83 	&	4.394E-02	&	1.86 	&	4.115E-02	&	1.88 	&	9.270E-02	&	1.86 	\\
0.8	&	   1/16 	&	2.716E-02	&	1.87 	&	1.168E-02	&	1.91 	&	1.087E-02	&	1.92 	&	2.502E-02	&	1.89 	\\
	&	   1/32 	&	7.161E-03	&	1.92 	&	3.020E-03	&	1.95 	&	2.804E-03	&	1.96 	&	6.549E-03	&	1.93 	\\
	&	   1/64 	&	1.843E-03	&	1.96 	&	7.689E-04	&	1.97 	&	7.125E-04	&	1.98 	&	1.679E-03	&	1.96 	\\
\hline
\end{tabular}
\end{center}
}
\end{table}
\subsection{Example of solution with weak regularity}

We apply the fractional $\theta$-methods to the Bagley-Torvik equation:
\begin{equation}\label{nume.4}
\begin{split}
\frac{\mathrm{d}^2 u}{\mathrm{d}x^2}+2 {}^C D_{0}^{\frac{3}{2}} u+2 u=f(x), \quad x \in (0,L],
\end{split}
\end{equation}
with initial conditions $u(0)=u'(0)=0$. ${}^C D_{0}^{\frac{3}{2}}$ denotes the Caputo differential operator, and under the initial conditions we have ${}^C D_{0}^{\frac{3}{2}}=\mathit{I}^{-\frac{3}{2}}$.
\par
We discretize the second order derivative term $\frac{\mathrm{d}^2 u}{\mathrm{d}x^2}$ by $\theta_1$-method and discretize ${}^C D_{0}^{\frac{3}{2}} u$ by $\theta_2$-method. Note that $\theta_1$ may be different from $\theta_2$. The exact solution is taken as $u(x)=x^{\mu}+x^5$ with $\mu>1$. Hence, $f(x)$ can be derived correspondingly. Let $L=1$.
For $\mu<4$, we use the formula (\ref{conv.12}) to obtain a second-order convergence rate. Refer \cite{Lubich1} for more information.
\par
In Table \ref{tab3} and Table \ref{tab4}, we take $\mu=1.1$ with different pairs $(\theta_1,\theta_2)$ for the fractional BT-$\theta$ method and the fractional BN-$\theta$ method, respectively.
The column Error$(\theta_1,\theta_2)$ represents the corresponding error derived by $\max_{1 \leq n \leq N}|U^n-u(x_n)|$, where $U^n$ is the numerical solution. The convergence rate is $O(h^2)$ in spite of different choices $(\theta_1,\theta_2)$.
\begin{table}[tbhp]
{\footnotesize
\caption{The convergence rate for the fractional BT-$\theta$ method with $\mu=1.1$.}\label{tab3}
\begin{center}
\begin{tabular}{|l|c|c|c|c|c|c|c|c|} \hline
$h$ & Err(0,0) & rate &  Err(-1,0.2) & rate &  Err(0.45,-0.1) & rate &  Err(-0.5,-2) & rate\\
\hline
   1/4  	&	5.184E-01	&	---     &	6.167E-01	&	---     &	3.908E-01	&	---     &	8.001E-01	&	--- \\
   1/8  	&	1.398E-01	&	1.89 	&	1.895E-01	&	1.70 	&	9.982E-02	&	1.97 	&	2.583E-01	&	1.63 	\\
   1/16 	&	3.812E-02	&	1.88 	&	5.714E-02	&	1.73 	&	2.670E-02	&	1.90 	&	8.206E-02	&	1.65 	\\
   1/32 	&	1.008E-02	&	1.92 	&	1.604E-02	&	1.83 	&	7.017E-03	&	1.93 	&	2.421E-02	&	1.76 	\\
   1/64 	&	2.598E-03	&	1.96 	&	4.266E-03	&	1.91 	&	1.805E-03	&	1.96 	&	6.663E-03	&	1.86 	\\
   1/128	&	6.600E-04	&	1.98 	&	1.101E-03	&	1.95 	&	4.584E-04	&	1.98 	&	1.754E-03	&	1.93 	\\
\hline
\end{tabular}
\end{center}
}
\end{table}
\begin{table}[tbhp]
{\footnotesize
\caption{The convergence rate for the fractional BN-$\theta$ method with $\mu=1.1$.}\label{tab4}
\begin{center}
\begin{tabular}{|l|c|c|c|c|c|c|c|c|} \hline
$h$ & Err(-0.2,-0.3) & rate &  Err(0,0) & rate &  Err(0.5,-0.1) & rate &  Err(1,0.7) & rate\\
\hline
   1/4  	&	7.201E-01	&	---     &	5.184E-01	&	---     &	5.897E-01	&	---     &	1.001E+00	&	---  \\
   1/8  	&	2.077E-01	&	1.79 	&	1.398E-01	&	1.89 	&	1.675E-01	&	1.82 	&	3.218E-01	&	1.64 	\\
   1/16 	&	5.893E-02	&	1.82 	&	3.812E-02	&	1.88 	&	4.727E-02	&	1.82 	&	9.887E-02	&	1.70 	\\
   1/32 	&	1.593E-02	&	1.89 	&	1.008E-02	&	1.92 	&	1.275E-02	&	1.89 	&	2.808E-02	&	1.82 	\\
   1/64 	&	4.153E-03	&	1.94 	&	2.598E-03	&	1.96 	&	3.320E-03	&	1.94 	&	7.522E-03	&	1.90 	\\
   1/128	&	1.061E-03	&	1.97 	&	6.600E-04	&	1.98 	&	8.477E-04	&	1.97 	&	1.948E-03	&	1.95 	\\
\hline
\end{tabular}
\end{center}
}
\end{table}
\section{Concluding remarks}\label{sec:conc}
In this paper, two families of novel fractional $\theta$-methods by constructing some new generating functions are proposed, the corresponding convergence, stability regions are developed, and some numerical tests are provided.
Specifically, the fractional BT-$\theta$ method connects FBDF2 with FTR while the fractional BN-$\theta$ methods links FBDF2 to GNGF2.
The convergence of the fractional BT-$\theta$ method is established directly by the linear multistep method while for the fractional BN-$\theta$ method, we derive the stability and consistency of the fractional convolution quadrature $\omega$ and then get the convergence of the method.
Both of the fractional $\theta$-methods result in a second-order convergence rate.
For an equation with a not regular solution, we can add some correction term to maintain the convergence rate.
We also discuss the stability regions in detail for the two fractional $\theta$-methods and illustrate the impact of different parameter $\theta$ on the stability regions.
Finally, numerical tests of our methods applied to the fractional ODE with smooth or nonsmooth solutions are implemented and the results confirm our theory.
\par
In another paper, we discuss some properties of the presented two families of novel fractional $\theta$ schemes and do some studies for fractional partial differential equations. In addition, we are considering other families of high order $\theta$ approximations.
\appendix
\section{Fast algorithm for the convolution weights}\label{sec:appendix}
We derive alternate formulas for the convolution weights $\omega_j$ in (\ref{thet.2}) and (\ref{thet.5}) which are obtained directly. Assume a generating function $\omega(\xi)$ takes the form
\begin{equation}\label{appe.1}
\omega(\xi)=p_1(\xi)^{\alpha}p_2(\xi)^{\beta}, \quad \alpha, \beta \in \mathbb{R},
\end{equation}
where $p_i$ ($i=1,2$) are polynomial functions with respect to $\xi$, and denote $p_{i,k}$ ($k=0,1,\cdots$) as the coefficients of $p_i(\xi)$.
\begin{theorem}\label{appe.the.1}
For the generating function $\omega(\xi)$ defined in (\ref{appe.1}),  the coefficients $\omega_k$ can be calculated recursively as follows
\begin{equation}\label{appe.2}
\omega_{k}=\frac{1}{k\psi_0}\sum_{j=0}^{k-1}(\phi_{k-j-1}-j\psi_{k-j})\omega_j,
\quad
\omega_0=p_{1,0}^{\alpha}p_{2,0}^{\beta},
\end{equation}
where $\phi_j$ and $\psi_j$ are the coefficients of $\phi(\xi):=\alpha p'_1(\xi) p_2(\xi)+\beta p_1(\xi)p'_2(\xi)$ and $\psi(\xi):=p_1(\xi)p_2(\xi)$, respectively.
Hence, $\phi_j=\alpha\sum_{s=0}^{j}p'_{1,j-s} p_{2,s}+\beta\sum_{s=0}^{j} p_{1,j-s}p'_{2,s}$, and $\psi_j=\sum_{s=0}^{j}p_{1,j-s}p_{2,s}$.
\end{theorem}
\textbf{Proof. } We take the first derivative of $\omega(\xi)$ and get
\begin{equation}\label{appe.3}
p_1 p_2 \omega'=(\alpha p'_1 p_2  +\beta p_1 p'_2) \omega.
\end{equation}
By the definitions of $\phi(\xi)$ and $\psi(\xi)$ and
considering both sides of equation (\ref{appe.3}) as functions, we expand with Taylor formula to obtain the $k$th coefficients, which satisfy
\begin{equation}\label{appe.4}
\sum_{j=0}^{k}\psi_{k-j}\omega'_k=\sum_{j=0}^{k}\phi_{k-j}\omega_k,
\end{equation}
where, $\sum_{j=0}^{\infty}\psi_{j}\xi^j :=\sum_{j=0}^{\infty}\frac{\psi^{(j)}(0)}{j !}\xi^j$,
$\sum_{j=0}^{\infty}\phi_{j}\xi^j :=\sum_{j=0}^{\infty}\frac{\phi^{(j)}(0)}{j !}\xi^j$ and
\begin{equation}\label{appe.5}
\sum_{j=0}^{\infty}\omega'_{j}\xi^j :=\sum_{j=0}^{\infty}\frac{\omega^{(j+1)}(0)}{j !}\xi^j
=\frac{\mathrm{d}}{\mathrm{d}\xi}\sum_{j=0}^{\infty}\frac{\omega^{(j)}(0)}{j !}\xi^{j}
=\sum_{j=1}^{\infty}j\omega_j\xi^{j-1}
=\sum_{j=0}^{\infty}(j+1)\omega_{j+1}\xi^{j}.
\end{equation}
Hence, $\omega'_{j}=(j+1)\omega_{j+1}$.
Now equation (\ref{appe.4}) can be formulated as
\begin{equation}\label{appe.6}
\sum_{j=0}^{k}(j+1)\psi_{k-j}\omega_{j+1}=\sum_{j=0}^{k}\phi_{k-j}w_j,
\end{equation}
and by extracting $\omega_{k+1}$, we have
\begin{equation}\label{appe.7}
\omega_{k+1}=\frac{1}{(k+1)\psi_0}\sum_{j=0}^{k}(\phi_{k-j}-j\psi_{k-j+1})\omega_j.
\end{equation}
The proof of the theorem is completed.
\begin{remark}
For the special case $p_2=1$, Weilbeer derived a corresponding formula in Theorem 5.3.1 in \cite{Weilbeer}. Note that for given polynomial functions $p_1$ and $p_2$, there are only finite many coefficients $p_{i,j}$ ($i=1,2$) that are non-zero. Hence, when computing weights $\{\omega_k\}_{k=0}^N$, the complexity of (\ref{appe.2}) is of $O(N)$ which is much more efficient than the direct calculation of (\ref{thet.2}) or (\ref{thet.5}).
\end{remark}
\begin{corollary}\label{appe.cor.1}
The convolution weights $\omega_k$ for the BT-$\theta$ method can be derived by the recursive formula
\begin{equation}\label{appe.8}\begin{split}
\omega_0=&\bigg(\frac{2-2\theta}{3-2\theta}\bigg)^{\alpha},
\quad
\omega_1=\frac{\phi_0\omega_0}{\psi_0},
\quad
\omega_2=\frac{1}{2\psi_0}[(\phi_0-\psi_1)\omega_1+\phi_1\omega_0],
\\
\omega_{k}=&\frac{1}{k\psi_0}\sum_{j=1}^{3}[\phi_{j-1}-(k-j)\psi_{j}]\omega_{k-j}, \quad k \geq 3,
\end{split}\end{equation}
where,
\begin{equation}\label{appe.9}
\phi_0=\frac{\alpha}{2}(2\theta^2-5\theta+4),
\quad \phi_1=\alpha(2\theta-1)(1-\theta),
\quad \phi_2=\frac{\alpha\theta}{2}(2\theta-1),
\end{equation}
and
\begin{equation}\label{appe.10}
\psi_0=\frac{1}{2}(3-2\theta)(1-\theta),
~\psi_1=\frac{1}{2}(1-2\theta)(3\theta-4),
~\psi_2=\frac{1}{2}(1-\theta)(1-6\theta),
~\psi_3=\frac{1}{2}\theta(1-2\theta).
\end{equation}
\end{corollary}
\begin{corollary}\label{appe.cor.2}
The convolution weights $\omega_k$ for the BN-$\theta$ method can be derived by the recursive formula
\begin{equation}\label{appe.11}\begin{split}
\omega_0=&\frac{2^{\alpha}(1-\alpha\theta)}{(3-2\theta)^{\alpha}},
\quad
\omega_1=\frac{\phi_0\omega_0}{\psi_0},
\quad
\omega_2=\frac{1}{2\psi_0}[(\phi_0-\psi_1)\omega_1+\phi_1\omega_0],
\\
\omega_{k}=&\frac{1}{k\psi_0}\sum_{j=1}^{3}[\phi_{j-1}-(k-j)\psi_{j}]\omega_{k-j}, \quad k \geq 3,
\end{split}\end{equation}
where,
\begin{equation}\label{appe.12}\begin{split}
\phi_0=&2\alpha(\theta-1)(\alpha\theta-1)-\alpha\theta(\theta-\frac{3}{2}),
\\
\phi_1=&\alpha(2\theta^2+3\alpha\theta-4\alpha\theta^2-1),
\\
\phi_2=&\alpha\theta(\frac{1}{2}-\theta-\alpha+2\alpha\theta),
\end{split}\end{equation}
and
\begin{equation}\label{appe.13}\begin{split}
\psi_0=&\frac{1}{2}(3-2\theta)(1-\alpha\theta),
\quad\psi_1=\frac{\alpha\theta}{2}(3-2\theta)+2(1-\theta)(\alpha\theta-1),
\\
\psi_2=&\frac{1}{2}(\alpha\theta-1)(2\theta-1)+2\alpha\theta(\theta-1),
\quad\psi_3=\frac{1}{2}\alpha\theta(1-2\theta).
\end{split}\end{equation}
\end{corollary}
\section*{Acknowledgments}
The authors are grateful to Professor Buyang Li for his valuable suggestions which improve the presentation of this work.


\end{document}